\documentclass{amsart}

\usepackage{amssymb,color}
\usepackage{amsfonts}
\usepackage{amsmath}
\usepackage{euscript}
\usepackage{enumerate}
\usepackage{pdfsync}
\newtheorem*{Theorem1'}{Theorem 1'}

\theoremstyle{definition}

\theoremstyle{remark}

\newcommand \Z{{\mathbb Z}}

\begin{document}

\title[The automorphism tower of the Mennicke group $M(-1,-1,-1)$]{The automorphism tower of the Mennicke group $M(-1,-1,-1)$}

\author{Fernando Szechtman}
\address{Department of Mathematics and Statistics, University of Regina, Canada}
\email{fernando.szechtman@gmail.com}
\thanks{The author was partially supported by NSERC grant RGPIN-2020-04062}





\begin{abstract} We compute the automorphism tower of the centerless Mennicke group $M(-1,-1,-1)$.
\end{abstract}

\maketitle

\section{Introduction}

In 1959, Mennicke \cite{Me} gave the first example of a finite group of deficiency zero requiring three generators.
A finite presentation $\langle X|R\rangle$ has deficiency $|X|-|R|$, and the deficiency of a finitely presented
group is the maximum of the deficiencies of its finite presentations. Mennicke's groups greatly contributed
to the interest in finite groups of deficiency zero, which can traced as far back as the dawn of the twentieth century to the work
of Miller \cite{Mi1, Mi2}. Following Mennicke's paper, many examples of finite groups of deficiency zero were found.
Among those requiring 3 generators are the groups discovered by Wamsley \cite{W2,W4},
Post \cite{P}, Johnson \cite{J}, Jamali \cite{Jam,Jam2}, and Allcock \cite{Al}. Amidst those needing 2 generators
are the groups found by Macdonald \cite{M}, Wamsley \cite{W3, W4}, Campbell and Robertson \cite{CR}, Robertson \cite{R},
Campbell, Robertson, and Thomas \cite{CRT}, Kenne \cite{K}, and Abdolzadeh and Sabzchi \cite{AS, AS2}.

Given integers $a,b,c$, the Mennicke group $M(a,b,c)=\langle x,y,z\,|\, x^y=x^a, y^{z}=y^b, z^x=z^c\rangle$.
Mennicke showed that $M(a,a,a)$ is finite provided $a\geq 2$. The abelianization of $M(a,b,c)$
makes it clear that $M(a,b,c)$ has deficiency zero and requires 3 generators whenever $a-1,b-1,c-1$ share a prime factor. The problem of the finiteness of the general Mennicke groups
$M(a,b,c)$ was studied by Macdonald and Wamsley \cite{W}, Schenkman~\cite{S}, and Jabara~\cite{Ja}. A sufficient condition
is that $a,b,c\not\in\{-1,1\}$. Upper bounds for the order of 
 $M(a,b,c)$ were provided by Johnson and Robertson \cite{JR}, Albar and Al-Shuaibi \cite{AA}, and Jabara \cite{Ja}.
The actual order of $M(a,b,c)$ is known only in certain cases studied in \cite{M,A,AA,Ja}.

In this paper, we investigate the automorphism tower of the infinite, centerless, and metabelian Mennicke group 
$$
M=M(-1,-1,-1)=\langle x,y,z\,|\, x^y=x^{-1}, y^{z}=y^{-1}, z^x=z^{-1}\rangle,
$$
as well as the related automorphism tower of a characteristic subgroup $V$ of $M$. We aim at studying
the structure and automorphism group of groups of deficiency zero (see \cite{MS,MS2}), and in the case
of $M$ we can go beyond and compute its full automorphism tower, which turns out to have length 2.
We note that the automorphism group of various types of metabelian groups have been
studied by several authors, see \cite{BBDM, BC, C, C2, CS, D, Di, G, GG, Ma, Ma2, Men}, for instance.
We recall, moreover, the celebrated result of Wielandt \cite{Wi} from 1939
to the effect that the automorphism tower of any finite centerless group terminates after finitely many steps.
We refer to Robinson's book \cite[Section 13.5.4]{Ro} for the definition of the automorphism tower of a centerless
group as well as a proof of Wielandt's theorem. His result was extended by Rae and Roseblade \cite{RR} to \v{C}ernikov groups in 1970, by Hulse \cite{H} to polycyclic groups in 1970,
and finally by Thomas \cite{T} in 1985, who showed that the automorphism tower of {\em any} centerless group eventually terminates
at some ordinal. It is possible to define the automorphism tower of any group, not necessarily centerless, and Hamkins \cite{Ha} proved in 1998
that the automorphism tower of {\em any} group actually terminates.

We begin by studying the structure of $M$. It turns out that $M$ is a centerless group such
that $[M,M]=M^2$ is a free abelian group of rank 3 and $M/[M,M]\cong C_2\times C_2\times C_2$ (in fact,
all factors of the lower central series of $M$ are isomorphic to $C_2\times C_2\times C_2$). In particular,
$M$ is metabelian and polycyclic. Moreover, all nontrivial elements of $M$ have infinite order, except for
those in the coset $xyzM^2$, of all whose elements have order 2. We use this information to determine $\mathrm{Aut}(M)$ and
$\mathrm{Aut}(\mathrm{Aut}(M))$, which is shown to be a complete group. 

The study of $\mathrm{Aut}(M)$ is carried out by means of the natural map $\Lambda:\mathrm{Aut}(M)\to 
\mathrm{Aut}(M/M^2)$. We show that 
$\{xM^2,yM^2,zM^2\}$ and $\{yxM^2,yzM^2,zxM^2\}$ are
$\mathrm{Aut}(M)$-orbits of $M/M^2$, where $M^2$ and $xyzM^2$ are clearly fixed.
Denoting by $\theta$ the automorphism of $M$ such that $x\mapsto y$, $y\mapsto z$, $z\mapsto x$,
we have that $\mathrm{Aut}(M)^\Lambda=\langle\theta\rangle^\Lambda$ is cyclic
of order 3. The kernel of $\Lambda$ is an extension of $\mathrm{Inn}(M)$ by $C_2\times C_2\times C_2$, and
$\mathrm{Out}(M)$ is the semidirect product of $C_2\times C_2\times C_2$ by $C_3$,
acting by cyclic permutation. It turns out that $G=\mathrm{Aut}(M)$ also has a characteristic subgroup $R$
that is free abelian of rank 3 and such that $G/R\cong C_2\times C_2\times C_2$, but $G\not\cong M$ as $G/G^2\cong C_2\times C_2$.

The study of $\mathrm{Aut}(G)$ is considerably more complicated. It transpires 
that $\mathrm{Out}(G)\cong C_2$ and that $\mathrm{Inn}(G)$ is a characteristic subgroup of $\mathrm{Aut}(G)$, so by a well-known result of Burnside, $\mathrm{Aut}(G)$ is a complete group.

The Mennicke group $M$ has 7 subgroups of index 2, but only one 
of these, say $V$, is 2-generated and torsion-free, so $V$ is characteristic
in $M$. The restriction map $\mathrm{Aut}(M)\to \mathrm{Aut}(V)$ is shown to be group monomorphism,
and we may view $\mathrm{Aut}(M)$ as a normal subgroup of $\mathrm{Aut}(V)$ with index~2 and trivial centralizer. 
This induces an imbedding $\mathrm{Aut}(V)\to
\mathrm{Aut}(\mathrm{Aut}(M))$, by conjugation, which is proven to be an isomorphism. Thus, the automorphism
tower of the centerless group $V$ terminates in the complete group $\mathrm{Aut}(V)$.

Presentations of all relevant groups, namely $V$, $\mathrm{Aut}(M)$, and $\mathrm{Aut}(\mathrm{Aut}(M))\cong \mathrm{Aut}(V)$
are given.

In terms of notation, given a group $T$ and $a,b\in T$, we set
$$
[a,b]=a^{-1}b^{-1}ab,\; b^a=a^{-1}ba,\; {}^a b=aba^{-1}.
$$
We write $T=\gamma_1(T),\gamma_2(T),\gamma_3(T),\dots$ for the terms of the lower central series of~$T$, so that 
$\gamma_{i+1}(T)=[T,\gamma_i(T)]$. Furthermore, for $n\in\Z$, we let $T^n$ stand for the subgroup of $T$ generated by
all $t^n$, with $t\in T$.  Function composition proceeds from left to right. If $\Lambda:S\to T$ is a group homomorphism,
we write $S^\Lambda=\{s^\Lambda\,|\, s\in S\}$, where $s^\Lambda$ indicates the image of $s$ under $\Lambda$.

\section{Finding $[M,M]$, $M^2$, and $M/[M,M]$}

Since $x^y=x^{-1}$, we have $x^{y^2}=x$, that is, $[x,y^2]=1$. Applying $\theta$ once yields $[y,z^2]=1$,
and a second application gives $[z,x^2]=1$. These 3 commutator identities imply $[x^2,y^2]=1$, $[y^2,z^2]=1$,
and $[z^2,x^2]=1$, whence $T=\langle x^2,y^2,z^2\rangle$ is an abelian group. We claim that 
$T$ is a normal subgroup of $M$. Indeed, from $(y^2)^x=y^2$ and $(z^2)^x=z^{-2}$ we deduce $T^x=T$.
As $T^\theta=T$, use of $\theta$ yields $T^y=T$ and $T^z=T$, which proves the claim.

The defining relations of $M$ imply
\begin{equation}
\label{sq}
[y,x]=x^2, [z,y]=y^2, [x,z]=z^2.
\end{equation}

One consequence of (\ref{sq}) is that $M/T$ is abelian, as $x,y,z$ commute modulo $T$. Thus, $[M,M]\subseteq T$.
But (\ref{sq}) also shows that $T\subseteq [M,M]$, so $T=[M,M]$. Thus, $M$ is metabelian.

Now, in general, if $G$ is any group, then $[G,G]\subseteq G^2$, and if $G=\langle X\rangle$, then $[G,G]=G^2$
if and only if $x^2\in [G,G]$ for every $x\in X$. Therefore, another consequence of (\ref{sq}) is that $M^2=[M,M]$. Let $\rho:M\to M/M^2$ be the canonical projection. As $M^\rho$ is abelian, generated by $x^\rho,y^\rho, z^\rho$,
and the square of each of these elements is trivial, we have $|M/M^2|\leq 8$, with equality if and only
if $M$ has an image of order 8 with $M^2$ in its kernel. One such image is $C_2\times C_2\times C_2$.
Thus $M/M^2\cong C_2\times C_2\times C_2$ and $M$ is disjoint union of the cosets $M^2, x M^2, yM^2, zM^2, xy M^2, yz M^2, zx M^2, 
xyz M^2$.

\section{A normal form for the elements of $M=\langle x \rangle\langle y \rangle\langle z \rangle$}

The defining relations of $M$ readily imply that $M=\langle x \rangle\langle y \rangle\langle z \rangle$.
We claim that every element of~$M$ can be written in exactly one way in the form $x^i y^j z^k$, where $i,j,k\in\Z$.
Indeed, we have group epimorphisms $f_1,f_2,f_3$ from $M$ onto
$D_\infty=\langle u,v\,|\, u^v=u^{-1}, v^2=1\rangle$, respectively given by
$$
x\mapsto u, y\mapsto v, z\mapsto 1;\; x\mapsto 1, y\mapsto u, z\mapsto v;\; x\mapsto v, y\mapsto 1, z\mapsto u.
$$
If $a,b,c\in\Z$ and $x^i y^j z^k=x^a y^b z^c$, then applying $f_1,f_2,f_3$
to both sides and using the uniqueness of expression in $\langle u,v\,|\, u^v=u^{-1}, v^2=1\rangle$ yields $i=a$, $j=b$, and $k=c$.
In particular, $M^2$ is a free abelian group with basis $\{x^2,y^2,z^2\}$. 

\section{Upper and lower central series of $M$}

We next claim that $Z(M)$ is trivial. Indeed, suppose $x^i y^j z^k\in Z(M)$. Applying $f_1,f_2,f_3$ above and using that
the center of $\langle u,v\,|\, u^v=u^{-1}, v^2=1\rangle$ is trivial, forces $i=j=k=0$.

Now $\gamma_2(M)=\langle x^2,y^2,z^2\rangle$. Suppose we have shown that 
$\gamma_i(M)=\langle x^{2^{i-1}},y^{2^{i-1}},z^{2^{i-1}}\rangle$ for some $i\geq 2$. Then $\gamma_{i+1}(M)$ is the normal
subgroup of $M$ generated by $[x,z^{2^{i-1}}]=z^{2^i}$, $[y,x^{2^{i-1}}]=x^{2^i}$, $[z,y^{2^{i-1}}]=y^{2^i}$. As
the subgroup of $M$ generated by these elements is already normal in $M$, we deduce 
$\gamma_{i+1}(M)=\langle x^{2^{i}},y^{2^{i}},z^{2^{i}}\rangle$. It follows that $\gamma_n(M)=
\langle x^{2^{n-1}},y^{2^{n-1}},z^{2^{n-1}}\rangle$ for all $n\geq 1$, with 
$$
\gamma_n(M)/\gamma_{n+1}(M)\cong \langle x^{2^{n-1}},y^{2^{n-1}},z^{2^{n-1}}\rangle/\langle x^{2^{n}},y^{2^{n}},z^{2^{n}}\rangle\cong
C_2\times C_2\times C_2.
$$

\section{The order of the elements of $M$}

If $u\in M$ is a nontrivial element of finite order $n$,
then $1=(u^n)^2=(u^2)^n$ with $u^2\in M^2$. As $M^2$ is free abelian, $u^2=1$ (so $n=2$). 
Now $u$ belongs to one and only one of the 8 cosets of $M^2$ described above, and it is certainly not in $M^2$.
Suppose, if possible, that $u\in x M^2$. Then $u=x v$, where $v=x^{2i} y^{2j} z^{2k}$ and $i,j,k\in\Z$. Thus
$$1=u^2=xvxv=x^2 v^x v=x^2 x^{2i} y^{2j} z^{-2k}x^{2i} y^{2j} z^{2k}=x^{2+4i} y^{4j},
$$
so $2+4i=0$, a contradiction. Thus $xM^2$ has no elements of finite order. Appealing to $\theta$ we see that
$yM^2$ and $zM^2$ have no elements of finite order either. Suppose next that $u\in xy M$. Then $u=xy v$, where $v=x^{2i} y^{2j} z^{2k}$ and $i,j,k\in\Z$.
Thus
$$1=u^2=xyvxyv=xyxy v^{xy} v = y^2 v^{xy} v=y^2 x^{-2i} y^{2j} z^{-2k}x^{2i} y^{2j} z^{2k}=y^{2+4j},
$$
so $2+4j=0$, a contradiction. Thus $xyM^2$ has no elements of finite order. Use of $\theta$ implies that
$yzM^2$ and $zxM^2$ have no elements of finite order either. Thus, the only possibility is that $u\in xyz M^2$, that is $u=xyz v$, where $v=x^{2i} y^{2j} z^{2k}$ and 
$i,j,k\in\Z$. Now
$$
u^2=xyz v xyz v=xyz xyz v^{xyz } v=(xyz)^2 x^{-2i} y^{-2j} z^{-2k}x^{2i} y^{2j} z^{2k}=(xyz)^2,
$$
where
$$
xyzxyz=xyxz[z,x]yz=x^2 y[y,x]z[z,x]yz=x^2 yz [y,x]^z yz [z,x]^{yz}=x^2 z^2 (x^2)^{zyz} (z^{-2})^{yz},
$$
so
$$
(xyz)^2=x^2 z^2 x^{-2} z^{-2}=1.
$$
Thus, the only nontrivial elements of $M$ of finite order are in $xyz M^2$, all of which have order 2.

\section{Conjugating formulas in $M$}

The following formulas will be used implicitly below. For $a,b,c\in\Z$, we have
$$
(z^c)^{x^a}=(z^{x^a})^c=(z^{(-1)^a})^c=z^{c (-1)^a}.
$$
Use of $\theta$ yields the formulas:
$$
(z^c)^{x^a}=z^{c (-1)^a}, (x^a)^{y^b}=x^{a (-1)^b}, (y^b)^{z^c}=y^{b (-1)^c}.
$$
As a consequence, we have
$$
(y^b)^{x^a}=x^{-a}y^b x^a=y^b y^{-b} x^{-a}y^b x^a=y^b x^{-a (-1)^b} x^a=y^b x^{a(1-(-1)^b)}.
$$
Application of $\theta$ then yields
$$
(y^b)^{x^a}=y^b x^{a(1-(-1)^b)}, (z^c)^{y^b}=z^c y^{b(1-(-1)^c)}, (x^a)^{z^c}=x^a z^{c(1-(-1)^a)}.
$$

\section{The stabilizers of $xM^2$, $yM^2$, and $zM^2$ in $\mathrm{Aut}(M)$ all act trivially on $M/M^2$}

As $M^2$ is a characteristic subgroup of $M$, we have an action of $\mathrm{Aut}(M)$ on $M/M^2$ by means of
automorphisms, that is, a group homomorphism $$\Lambda:\mathrm{Aut}(M)\to 
\mathrm{Aut}(M/M^2)\hookrightarrow\mathrm{Sym}(M/M^2).$$

If $\alpha\in \mathrm{Aut}(M)$ and $u\in M$, we say that $\alpha$ fixes $uM^2$ if $(uM^2)^\alpha=uM^2$.
Every automorphism of $M$ fixes $xyzM^2$, as it is the only coset of $M^2$ all of whose elements have order 2.

Let $K=\ker(\Lambda)$. Clearly $K$ is contained in the stabilizer of $xM^2$. We claim that, in fact, they are equal to each other, that is, the stabilizer of $xM^2$ actually fixes $yM^2$ and $zM^2$ as well.

Let $\alpha$ be an automorphism of $M$ that fixes $xM^2$, and set $x_0=x^\alpha$, $y_0=y^\alpha$. Then $x_0=x u$, 
where $u\in M^2$, so that $u=x^{2a}y^{2b}z^{2c}$ and
$$
x_0^{y_0}=x_0^{-1}=xx^{-2-2a}y^{-2b}z^{2c}.
$$
Clearly, $\alpha$ sends $yM^2$ to one of the following cosets: $yM^2,zM^2,xyM^2,yzM^2,zxM^2$. 
Suppose, if possible, that $yM^2$ is sent to $zM^2$. Then $y_0=z  v$, where $v\in M^2$, so that $v=x^{2d}y^{2e}z^{2f}$,
and therefore
$$
x_0^{y_0}=xx^{2a}y^{-2b}z^{2+2c+4f}.
$$
Since $M^2$ is free with basis $x^2,y^2,z^2$, we deduce $2a=-2-2a$, which is impossible.

Suppose next that $yM^2$ is sent to $zxM^2$. Then $y_0=zx  v$, where $v\in M^2$, so that $v=x^{2d}y^{2e}z^{2f}$,
whence
$$
x_0^{y_0}=xx^{2a}y^{-2b}z^{-2+2c+4f}.
$$
We see, as above, that this is absurd.

Suppose next that $yM^2$ is sent to $yzM^2$. Then $y_0=yzv$, where $v\in M^2$, so that $v=x^{2d}y^{2e}z^{2f}$,
which implies
$$
x_0^{y_0}=xx^{-2-2a}y^{-2b}z^{2+2c+4f}.
$$
Thus, $2+4f=0$. Since $f$ is an integer, this is impossible.

Thus $\alpha$ sends $yM^2$ to  $yM^2$  or $xyM^2$. Since $\alpha$ induces an automorphism of $M/M^2$ that
fixes $xM$ and $xyzM$, it follows
that either $\alpha$ fixes both $yM^2$ and $zM^2$, or else sends $yM^2$ to $xyM^2$ and $zM^2$ to $zxM^2$.

Suppose, if possible, that $\alpha$ sends $zM^2$ into $xzM^2$. Then $z_0=z^\alpha= xz v$, where $v\in M^2$, so that 
$v=x^{2d}y^{2e}z^{2f}$. We must have $z_0^{x_0}=z_0^{-1}$. But
$$
z_0^{x_0}=xz x^{2d}y^{2e+4b}z^{-2-2f+4c}, z_0^{-1}=xz x^{2-2d}y^{2e}z^{2f}.
$$
Therefore $2d=2-2d$, which is absurd.

Thus $\alpha$ must fix $yM^2$ and $zM^2$. This proves the claim. Since $\mathrm{Stab}(xM^2)=K$, use of $\theta$ yields $\mathrm{Stab}(y M^2)= K=\mathrm{Stab}(zM^2)$ as well.

\section{The $\mathrm{Aut}(M)$-orbits of $M/M^2$}

Since every automorphism of $M$ fixes $M^2$ and $xyz M^2$, it follows that
$\mathrm{Aut}(M)$ permutes the remaining 6 cosets of $M^2$. 

Under the permutation action of
$\mathrm{Aut}(M)$ on $M/M^2$, use of $\theta$ reveals that $xM^2,yM^2,zM^2$ are in the same orbit
and so are $xyM^2,yz M^2,zxM^2$. Thus, either these are $\mathrm{Aut}(M)$-orbits, or so is
their union. We claim that $\{xM^2,yM^2,zM^2\}$ and  $\{xyM^2,yz M^2,zxM^2\}$ are $\mathrm{Aut}(M)$-orbits.

Let $\phi\in \mathrm{Aut}(M)$ be arbitrary. It is not possible for $\phi$ to send all of $\{xM^2,yM^2,zM^2\}$ into 
$\{xyM^2,yz M^2,zxM^2\}$,
because $\phi$ induces an automorphism of $M/M^2$, so the images of the generators $xM^2,yM^2,zM^2$ of $M/M^2$ under $\phi$
must also generate $M/M^2$, and $xyM^2,yz M^2,zxM^2$ do not generate $M/M^2$. Thus, $\phi$ must send at least
one of $xM^2,yM^2,zM^2$ back into $\{xM^2,yM^2,zM^2\}$. But then $\phi\theta^i$ fixes one of $xM^2,yM^2,zM^2$
for some $i\in\Z$, so by above $\phi\theta^i$ fixes all of $xM^2,yM^2,zM^2$, whence $\phi$ stabilizes
$\{xM^2,yM^2,zM^2\}$, and hence $\{xyM^2,yz M^2,zxM^2\}$, as claimed.

\section{The image and kernel of $\Lambda:\mathrm{Aut}(M)\to\mathrm{Aut}(M/M^2)$}

We claim that $\mathrm{Aut}(M)^\Lambda=\langle \theta\rangle^\Lambda$ is a cyclic group of order 3.
Indeed, by above, $\mathrm{Aut}(M)$ permutes the elements $xM,yM,zM$ of $M/M^2$. It 
permutes them transitively by means of $\langle\theta\rangle$, so 
$$\mathrm{Aut}(M)=\mathrm{Stab}(xM^2) \langle\theta\rangle= K \langle\theta\rangle,
$$
whence $\mathrm{Aut}(M)^\Lambda=\langle \theta\rangle^\Lambda$ is a cyclic group of order 3. A consequence to be used later
is that the permutation action of $\mathrm{Aut}(M)$ on $\{xy M^2, yz M^2, zx M^2\}$ is given by a 3-cycle. In particular,
no automorphism of $M$ induces a transposition on $\{xy M^2, yz M^2, zx M^2\}$.

Let us determine the kernel $K$ of $\Lambda$. By definition, $K$ consists of all automorphisms of $M$ that act like the identity
on $M/M^2$.  As $M/M^2$ is abelian, the group $\mathrm{Inn}(M)$ acts trivially on $M/M^2$, so $\mathrm{Inn}(M)$ 
is included in $K$. But $K$ strictly contains $\mathrm{Inn}(M)$.
Indeed, let $\phi\in K$. Then $x_0=x^\phi=x x^{2a}y^{2b}z^{2c}$ and $y_0=y^\phi=y x^{2d}y^{2e}z^{2f}$.
We must have $x_0^{y_0}=x_0^{-1}$. But
$$
x_0^{y_0}=x^{-1} x^{-2a}y^{2b}z^{2c+4f}, x_0^{-1}=x^{-1}x^{-2a} y^{-2b}z^{2c},
$$
which implies that $b=0=f$. A similar argument involving $y^\theta$ and $z^\theta$ shows that
$$
x^\phi=x x^{2a}z^{2c}, y^\phi=y x^{2d}y^{2e}, z^\phi=z y^{2g}z^{2h}.
$$
But then 
$$
(x^2)^\phi=(x^2)^{1+2a}, (y^2)^\phi=(y^2)^{1+2e}, (z^2)^\phi=(z^2)^{1+2h}.
$$
As $\phi$ induces an automorphism of $M^2$, which has basis $x^2,y^2,z^2$, it follows that $1+2a,1+2e,1+2h$ are units in $\Z$,
so each of them is equal to $\pm 1$. Thus $a,e,h\in\{0,-1\}$.

Thus
\begin{equation}\label{p0}
x^\phi=x^{\pm 1}z^{2c}, y^\phi=y^{\pm 1} x^{2d}, z^\phi=z^{\pm 1} y^{2g}.
\end{equation}
Following $\phi$ by conjugation by $x$, $y$, or $z$, if necessary, we may assume loss that
\begin{equation}\label{p}
x^\phi=x z^{2c}, y^\phi=y x^{2d}, z^\phi=z y^{2g}.
\end{equation}
These assignment are easily seen to preserve the defining relations of $M$, so they do extend to endomorphisms of $M$.

For $(c,d,g)\in\Z^3$, consider the endomorphism $P_{(c,d,g)}$ of $M$ defined by (\ref{p}). 
Addition of triples corresponds to composition of endomorphisms, and $(0,0,0)$ corresponds to the identity,
so each assignment (\ref{p}) extends to an automorphism of $M$, and hence so do all assignments (\ref{p0}). 

What is $K/\mathrm{Inn}(M)$? Given an arbitrary $\phi\in K$ using suitable inner automorphisms
we may assume that $\phi$ is given by (\ref{p}), as indicated above.
Further conjugating by suitable elements of $M^2$, we may use the above conjugating formulas of $M$
to subtract any multiple of 4 from $2c,2d,2g$. Thus, any $\phi\in K$ is congruent modulo $\mathrm{Inn}(M)$ to one of the form
$$
x\mapsto x z^{r}, y\mapsto y x^{s}, z\mapsto z y^{t},
$$
where $r,s,t\in\{0,2\}$. The above conjugating formulas in $M$ reveal that none of these automorphisms are inner, except when $r=s=t=0$.
They commute and their squares are all inner automorphisms, so $K/\mathrm{Inn}(M)\cong C_2\times C_2\times C_2$.

\section{Description of $\mathrm{Aut}(M)$ and $\mathrm{Out}(M)$}

By above, $\mathrm{Aut}(M)=K\rtimes\langle\theta\rangle$. Here $\mathrm{Inn}(M)\subset K$, so
$$
\mathrm{Out}(M)=\mathrm{Aut}(M)/\mathrm{Inn}(M)\cong (K/\mathrm{Inn}(M))\rtimes C_3\cong (C_2\times C_2\times C_2)
\rtimes C_3,
$$
where $C_3$ acts on $C_2\times C_2\times C_2$ by cyclic permutation. Here $M\cong \mathrm{Inn}(M)$ as $Z(M)$ is trivial. 

Let $X,Y,Z$ be the inner automorphisms corresponding to $x,y,z$, let $A,B,C$ be respectively defined by
$$
x\mapsto x z^{2}, y\mapsto y, z\mapsto z; x\mapsto x, y\mapsto y x^{2}, y\mapsto y; x\mapsto x, y\mapsto y, z\mapsto z y^{2},
$$
and let $D=\theta$. Then every element of $\mathrm{Aut}(M)$ can be written in one
and only one way in the form
$$
X^a Y^b Z^c A^i B^j C^k D^\ell,\quad a,b,c\in\Z, 0\leq i,j,k\leq 1, 0\leq\ell\leq 2.
$$

Note that $\mathrm{Aut}(M)$ is generated by $X,A,D$ but we prefer to give a presentation in terms of
the generators $X,Y,Z,A,B,C,D$, with defining relations:
\begin{equation}\label{defG1}
X^Y=X^{-1}, X^D=Y, Y^D=Z, Z^D=X, D^3=1, 
\end{equation}
\begin{equation}\label{defG2}
[A,B]=1, A^D=B, B^D=C, C^D=A, A^2=Z^2, X^A=XZ^2, Y^A=Y, Z^A=Z.
\end{equation}
The following relations are consequence of the previous ones:
$$
Y^Z=Y^{-1}, Z^X=Z^{-1}, [B,C]=1, [C,A]=1, B^2=X^2, C^2=Y^2, X^B=X, 
$$
$$
Y^B=YX^2, Z^B=Z, X^C=X, Y^C=Y, Z^C=ZY^2.
$$
Note that $\langle A,B,C\rangle$ is a free abelian group of rank 3. 

\section{The characteristic subgroup $V$ of $M$}

Let 
$$V=\langle xy,yz,zx,x^2,y^2,z^2\rangle=\langle xy,yz,zx\rangle=\langle xy,yz\rangle=\langle yz,zx\rangle=\langle zx,xy\rangle.
$$
Since  $V$ contains $[M,M]$, it is a normal subgroup of $M$. Clearly $[M:V]=2$.
There are 7 subgroups of $M$ of index 2, all necessarily containing $M^2$. They correspond to the 7 hyperplanes
of the $F_2$-vector space $M/M^2$. Three of these subgroups contain $xyz$, namely $\langle x,yz\rangle M^2$,
$\langle y,zx\rangle M^2$, and $\langle z,xy\rangle M^2$. They  form  an orbit under $\langle D\rangle$.
Another such orbit is formed by $\langle x,y\rangle M^2$,
$\langle x,z\rangle M^2$, and $\langle y,z\rangle M^2$. The seventh such subgroup is $V$, which is clearly
$\langle D\rangle$-invariant, and is thus the only one that could be $\mathrm{Aut}(M)$-invariant.
But $V$ is certainly $X$-invariant (being normal) and $D$-invariant (as just mentioned), and we readily verify that $V$ is $A$-invariant. Thus $V$ is in fact $\mathrm{Aut}(M)$-invariant.

Let us see that $V$ is a characteristic subgroup of $M$ by means structural properties, without appealing to
our knowledge of $\mathrm{Aut}(M)$. Now $\mathrm{Aut}(M)$ permutes the 7 subgroups of $M$ of index~2.
It cannot map $V$ into any of $\langle x,yz\rangle M^2$,
$\langle y,zx\rangle M^2$, and $\langle z,xy\rangle M^2$, as $V$ is torsion free and these 3 subgroups are not.
As indicated above $V$ can be generated by 2 elements. The subgroups $\langle x,y\rangle M^2$,
$\langle x,z\rangle M^2$, and $\langle y,z\rangle M^2$ form a $\langle D\rangle$-orbit, so if one of
them requires 3 generators, they all do. But $W=\langle x,y\rangle M^2=\langle x,y,z^2\rangle$
requires 3 generators, as $W^2=\langle x^2,y^2,z^4\rangle$ and $W/W^2\cong C_2\times C_2\times C_2$
cannot be generated by 2 elements.

We note that $V$ is generated by $u=xy$, $v=yz$, and $w=zx$ of $V$, subject to the defining relations:
$$
[u,v]=w^2 u^{-2} v^2, [v,w]=u^2 v^{-2} w^2, [w,u]=v^2 w^{-2} u^2,
$$
$$
uvw=w^2u^2v^{-2}, vwu=u^2v^2w^{-2}, wuv=v^2w^2u^{-2},
$$
$$
(u^2)^v=u^{-2}, (u^2)^w=u^{-2}, (v^2)^w=v^{-2}, (v^2)^u=v^{-2}, (w^2)^u=w^{-2}, (w^2)^v=w^{-2}.
$$
The following relations are consequence of the above:
$$
[u^2,v^2]=1, [v^2,w^2]=1, [w^2,u^2]=1.
$$

Observe that $V$ is a torsion-free group, as the nontrivial torsion elements of $M$ are all in $xyz M^2$ and none of these are in $V$.
As $V$ is a subgroup of $M$, it is also metabelian. Note as well that $V$ is centerless.
Indeed, consider the homomorphism $M\to \langle a,b\,|\, a^b=a^{-1}, b^2=1\rangle$, given by
$$
x\mapsto a, y\mapsto b, z\mapsto 1.
$$
Then the restriction to $V$, say $f$, satisfies
$$
xy\mapsto ab, yz\mapsto b, zx\mapsto a.
$$
It follows that $Z(V)\subseteq \ker(f)$. Now $V$ is the disjoint union of $V^2, xy V^2, yz V^2, zx V^2$, and applying
$f$ to elements of these cosets, we see that $\ker(f)$ is contained in $V^2$, and in fact in $\{y^{2j}z^{2k}\,|\, j,k\in\Z\}$,
as $a^i b^j=1$ if and only if $i=0$ and $j$ is even. Conjugating the central element
$y^{2j}z^{2k}$ by $xy$ we find that $k=0$, and conjugating the central element
$y^{2j}$ by $yz$ we deduce that $j=0$, so $Z(V)$ is trivial.

\section{The restriction map $\mathrm{Aut}(M)\to \mathrm{Aut}(V)$ is a group monomorphism}

Taking into account that $V$ is a characteristic subgroup of $M$, we may now consider the restriction map 
$\Gamma:\mathrm{Aut}(M)\to \mathrm{Aut}(V)$. 

We claim that $\Gamma$ is injective. Indeed, let $\psi\in\ker(\Gamma)$. Then $\psi$ acts like the identity on $V$ and hence on 
$M^2$. Now
$\psi=\phi D^\ell A^i B^j C^k$, with $\phi\in \langle X,Y,Z\rangle$, as indicated above. 
Since $A,B,C$ already act like the identity on $M^2$,
it follows that $\phi D^\ell$ acts like the identity of $M^2$.  Writing $\phi=X^a Y^b Z^c$ it follows easily
that $a,b,c$ must be even and $\ell=0$. We conclude that the kernel of the restriction map 
$$
\mathrm{Aut}(M)\to \mathrm{Aut}(M^2)
$$
is equal to the abelian group $\langle X^2,Y^2,Z^2, A,B,C\rangle$. Going back to $\psi\in\ker(\Gamma)$,
we have that $\psi= X^{2a} Y^{2b} Z^{2c} A^i B^j C^k$. This automorphism of $M$ must fix $xy$.
But
$$
(xy)^{X^{2a} Y^{2b} Z^{2c} A^i B^j C^k}=xy x^{2j+4a}z^{2i+4c}.
$$
Using the uniqueness of expression in $M$, we see that $j+2a=0$ and $i+2c=0$.
But $0\leq i,j\leq 1$, so $i=j=a=c=0$. Using that $\psi$ fixes $yz$ and $zx$ we deduce that $\psi=1$.

\section{The image and kernel of $\Pi:\mathrm{Aut}(V)\to \mathrm{Aut}(V/V^2)$}

Notice that $\langle x^2,y^2,z^2\rangle\subseteq V^2\subseteq M^2=\langle x^2,y^2,z^2\rangle$, so $V^2=M^2$,
and $V/V^2\cong C_2\times C_2$. To find $\mathrm{Aut}(V)$, 
we look at the natural map $\Pi:\mathrm{Aut}(V)\to \mathrm{Aut}(V/V^2)\cong S_3$. 
We claim that the image of $\Pi$ is the full automorphism group of $V/V^2$. Indeed, first of all, viewing $\mathrm{Aut}(M)$
as a subgroup of $\mathrm{Aut}(V)$, it follows from above that $\mathrm{Aut}(M)^\Pi$ contains a 3-cycle but no transpositions. To
see that $\mathrm{Aut}(V)^\Pi$ contains a transposition, consider the assignment
\begin{equation}\label{pis}
u\mapsto u w^2, v\mapsto w v^2, w\mapsto vu^2.
\end{equation}
The defining relations of $V$ allow us to verify that this extends to an endomorphism, say $\Psi$, of $V$ such that
$$
u^2\mapsto u^2, v^2\mapsto w^2, w^2\mapsto v^2,
$$
$$
V^2\mapsto V^2, uV^2\mapsto uV^2, vV^2\mapsto wV^2, wV^2\mapsto vV^2.
$$
It follows that $\Psi$ is an automorphism, as required. The inverse $\Psi^{-1}$ is given by
$$
u\mapsto u v^{-2}, v\mapsto w u^{-2}, w\mapsto vw^{-2}.
$$

We next claim that the kernel of $\Pi:\mathrm{Aut}(V)\to \mathrm{Aut}(V/V^2)$ consists of the restriction to $V$
of automorphisms of $M$ in the kernel of $\mathrm{Aut}(M)\to \mathrm{Aut}(M/M^2)$. Indeed, let $\phi$ be in the kernel of $\Pi$. Set
$$
u=xy, v=yz, w=zx, u_0=(xy)^\phi, v_0=(yz)^\phi, w_0=(zx)^\phi.
$$
Then
$$
u_0=xy x^{2i}y^{2j}z^{2k}, v_0=yz x^{2a}y^{2b}z^{2c},   w_0=zx x^{2d}y^{2e}z^{2f}.
$$
The relation $u^v=u w^2 u^{-2}v^2$ must be preserved by $\phi$, so we must have $u_0^{v_0}=u_0 w_0^2 u_0^{-2}v_0^2$.
One the one hand, 
$$
u_0^2=y^{2(1+2j)}, v_0^2=z^{2(1+2c)}, w_0^2=x^{2(1+2d)},
$$
so that
$$
u_0 w_0^2 u_0^{-2}v_0^2=xy x^{2i}y^{2j}z^{2k} x^{2(1+2d)}y^{-2(1+2j)}z^{2(1+2c)},
$$
while on the other hand,
$$
u_0^{v_0}=xy x^2 y^{-2} z^{4c} z^{2} x^{4a} x^{-2i} y^{-2j} z^{2k}.
$$
Looking at the exponents of $x^2$, we deduce that 
$$
i+d-a=0.
$$
Likewise, the relation $v^w=v u^2 v^{-2}w^2$ must be preserved by $\phi$. On the one hand, we have
$$
v_0 u_0^2 v_0^{-2}w_0^2=yz x^{2a}y^{2b}z^{2c} x^{2(1+2d)}y^{2(1+2j)}z^{-2(1+2c)},
$$
while on the other hand,
$$
v_0^{w_0}=yz y^2 z^{-2} x^{4d} x^{2} y^{4e} y^{-2b} z^{-2c} x^{2a}.
$$
Considering the exponents of $y^2$, we infer that 
$$
b+j-e=0.
$$
Finally, the relation $w^u=w v^2 w^{-2}u^2$ must be preserved by $\phi$. On the one hand, we have
$$
w_0 v_0^2 w_0^{-2}u_0^2=zx x^{2d}y^{2e}z^{2f} x^{-2(1+2d)}y^{2(1+2j)}z^{2(1+2c)},
$$
while on the other hand,
$$
w_0^{u_0}=zx z^2 x^{-2} y^{4j} y^{2} z^{4k} z^{-2f} x^{-2d} y^{2e}.
$$
An examination of the exponents of $z^2$ reveals that 
$$
f+c-k=0.
$$
But $\phi$ induces an automorphism on the free abelian group $V^2$ with basis $\{x^2,y^2,z^2\}$, and
$$
u_0^2=y^{2(1+2j)}, v_0^2=z^{2(1+2c)}, w_0^2=x^{2(1+2d)}.
$$
This implies that $j,c,d\in\{0,-1\}$. Now if $j=-1$, we may follow $\phi$ by conjugation by $z$  (restricted to $V$)
and make the new $j=0$, while keeping $c,d$ the same.  Likewise, if $c=-1$, we may follow $\phi$ by conjugation by $x$  
(restricted to $V$)
and make the new $c=0$, while keeping $j,d$ the same. Finally, if $d=-1$, we may follow $\phi$ by conjugation by $y$  (restricted to $V$)
and make the new $d=0$, while keeping $c,j$ the same. Thus, we may assume without loss that $j=c=d=0$, in which case
$$
i=a, b=e, f=k.
$$
This kind of automorphism is obtained by restricting to $V$ the
of automorphisms of $M$ in the kernel of $\mathrm{Aut}(M)\to \mathrm{Aut}(M/M^2)$ found above. This proves the claim.

We deduce from above that we may view $\mathrm{Aut}(M)$ as a subgroup of index 2, and hence normal, in $\mathrm{Aut}(V)$.

\section{$\mathrm{Inn}(M)$ is not a characteristic subgroup of $\mathrm{Aut}(M)$}

We claim that $\mathrm{Inn}(M)$ is not invariant under the automorphism of $\mathrm{Aut}(M)$ given by
conjugation by $\Psi$ (as defined in (\ref{pis})). Suppose, on the contrary, that this is false. For $t\in M$, let 
$C_t$ stand for conjugation by $t\in M$ {\em restricted to $V$}. Then, given any $t\in M$ 
there is a unique $t'\in M$ such that
$C_t^\Psi=C_{t'}$. The map $t\mapsto t'$ is an automorphism, say $\Delta$, of $M$, so that
$$
C_t^\Psi=C_{t^{\Delta}},\quad t\in M,
$$
and therefore
$$
C_{t^\Psi}=C_t^\Psi=C_{t^{\Delta}},\quad t\in V.
$$
As $V$ is centerless, this implies that $\Psi$ extends to the automorphism $\Delta$ of $M$, which is impossible
as no automorphism of $M$ induces a transposition on $\{xy M^2, yz M^2, zx M^2\}$.

It follows that $\mathrm{Inn}(M)$ is not a characteristic subgroup of $\mathrm{Aut}(M)$.

\section{The centralizer of $\mathrm{Aut}(M)$ in $\mathrm{Aut}(V)$ is trivial}

We claim that centralizer, say $Q$, of $\mathrm{Aut}(M)$ in $\mathrm{Aut}(V)$ is trivial. Indeed, as $Z(M)$ is trivial
we infer that $\mathrm{Aut}(M)$ is centerless. Suppose, if possible, that $Q$ is nontrivial. As $Q\cap \mathrm{Aut}(M)$
is trivial and $\mathrm{Aut}(V)$ is the disjoint union of $\mathrm{Aut}(M)$ and $\mathrm{Aut}(M)\Psi$, we deduce
that $\Upsilon\Psi\in Q$ for some $\Upsilon\in\mathrm{Aut}(M)$. This means that conjugation by $\Psi$ and $\Upsilon^{-1}$
agree on $\mathrm{Aut}(M)$. In particular, $\mathrm{Inn}(M)$ is invariant under conjugation by $\Psi$, a contradiction.

\section{The isomorphism $\mathrm{Aut}(V)\to \mathrm{Aut}(\mathrm{Aut}(M))$}

By above, we have a group monomorphism $\Omega:\mathrm{Aut}(V)\to \mathrm{Aut}(\mathrm{Aut}(M))$ by conjugation.
In addition, we know that $[\mathrm{Aut}(V):\mathrm{Aut}(M)]=2$. As shown below, the index of $\mathrm{Aut}(M)$ in
$\mathrm{Aut}(\mathrm{Aut}(M))$ is also 2, so $\Omega$ is an isomorphism.

\section{The characteristic subgroup $R$ of $\mathrm{Aut}(M)$}

Setting $G=\mathrm{Aut}(M)$, we proceed to determine $\mathrm{Aut}(G)$. 
We may view $M$ as a subgroup of $G$ via $x\mapsto X,y\mapsto Y,z\mapsto Z$ (that is,
via the natural imbedding $M\to\mathrm{Inn}(M)$). Routine calculations show that
$$
[G,G]=\langle XY,YZ,ZX,AB,BC,CA\rangle,\; G/[G,G]\cong C_2\times C_2\times C_3,
$$
$$
G^2=\langle XY,YZ,ZX,AB,BC,CA,D\rangle, G/G^2\cong C_2\times C_2,
$$
$$
[G,G]^2=\langle X^2,Y^2,Z^2\rangle=M^2, G/[G,G]^2\cong C_2^6\rtimes C_3,
$$
where $C_3$ acts on $C_2^6$ via two 3-cycles. Thus, there is a unique (normal) Sylow 2-subgroup, say $S$, of 
$G/[G,G]^2$, and $M/M^2$ is inside it. Here 
$$
S=U/[G,G]^2, U=\langle X,Y,Z,A,B,C\rangle,
$$
where $U$ is a characteristic subgroup of $G$, being the only subgroup of $G$ containing $[G,G]^2$ and having index 3 in $G$.
Note that $S$ is a 6-dimensional vector space over $F_2$.

We claim that $R=\langle A,B,C\rangle$ is the only normal abelian subgroup of $G$ containing $M^2$ 
and such that $R/M^2$ is a 3-dimensional
subspace of $S$. This makes $R$ into a characteristic subgroup of $G$. To prove the claim,
let $Q$ be any normal abelian subgroup of $G$ containing $M^2$ 
such that $R/M^2$ is a 3-dimensional
subspace of $S$. Then $Q$ contains elements of the form $u=X^aY^bZ^cA^dB^eC^f$, each exponent being 0 or 1,
the sum of which we call the length of $u$. If these elements contain terms only from $\{A,B,C\}$ or only from $\{X,Y,Z\}$, then $Q=R$, as $M$ is not abelian.
Suppose, if possible, that $Q$ contains elements with mixed terms, and let $m>1$ be the smallest 
length among these elements. We divide the analysis
into five cases:

$\bullet m=2$. As $Q$ is normal, then $XA,YB\in Q$, or $XB,YC\in Q$, or $XC,YA\in Q$, all of which contradict the fact that 
$Q$ is abelian.

$\bullet m=3$. As $Q$ is normal, then $XAB,YBC\in Q$, or $XBC,YAC\in Q$, or $XAC,YAB\in Q$, 
or $XYA,YZB\in Q$, or $XYB,YZC\in Q$, or $XYC,YZA\in Q$, all of which contradict the fact that 
$Q$ is abelian.

$\bullet m=4$. As $Q$ is normal, then $XYZA,XYZB, XYZC\in Q$, or $XABC,YABC, ZABC\in Q$, or
$Q$ contains an element $u$ with two factors from $\{X,Y,Z\}$ and two factors
from $\{A,B,C\}$. In the first case, $Q$ contains $A,B,C, XYZ$ and is at least 4-dimensional.
In the second case, $Q$ contains $X,Y,Z, ABC$ and is at least 4-dimensional. In the third case, $Q$
contains $uu^D$, which yields an element of $Q$ with mixed terms and of length 2. All of these cases
are impossible.

$\bullet m=5$. As $Q$ is normal, then either $XYABC,YZABC,XZABC\in Q$, in which case 
$Q$ contains the centeress group $V$, or else $XYZAB,XYZBC,XYZAC\in Q$, in which case $Q$ contains the non-commuting
elements $XYZAB, AB$. Both cases are absurd.

$\bullet m=6$. This is impossible, because then $Q/M^2$ is one dimensional.

This proves that $R$ is a characteristic subgroup of $G$. 

\section{The $\mathrm{Aut}(G)$-orbit of $M$}

We know that $M$ is not a characteristic subgroup of $G$, and we single out the automorphism, say $E$, of $G$ defined by
$A,B,C,D$ and satisfying
$$
A\mapsto A, B\mapsto B, C\mapsto C, D\mapsto D, X\mapsto XA,Y\mapsto YB, Z\mapsto ZC.
$$
Note that $E^2$ is conjugation by $ABC$. 

We claim that the $\mathrm{Aut}(G)$-orbit of $M$ consists of $M$ and $M^E$ only.
To prove the claim, let $Q$ be any normal subgroup of $G$ such that $[Q,Q]=M^2$ and $Q/M^2$ is a 3-dimensional subspace of~$S$,
and suppose that $Q$ is in the $\mathrm{Aut}(G)$-orbit of $M$. Let $\alpha$ be any automorphism of $G$ such that $M^\alpha=Q$,
and set $X_0=X^\alpha$ and $D_0=D^\alpha$. 
Now $U=MR$ and as $U$ and $R$ are characteristic subgroups of $G$, we infer $U=QR$. Here $U/R$ is a 3-dimensional
vector space over  $F_2$, and $\alpha$ induces
a vector space automorphism of $U/R$ that sends the basis $\{X R, X^D R, X^{D^2} R\}$ into the basis
$\{X_0 R, X_0^{D_0} R, X_0^{D_0^2} R\}$. It follows that $X_0$ cannot be congruent modulo $R$ to any of 
$XY,YZ,ZX, XYZ$, so these two bases are identical (we are, of course, disregarding the order of the basis vectors).
Thus, replacing $\alpha$ by an inner automorphism of $G$ associated to $\langle D\rangle$, we may assume
that $X_0\equiv X\mod R$. As $Q$ contains $M^2$, it follows that one of the following elements is in $Q$:
$X,XA,XB,XC,XAB,XAC,XBC, XABC$. In the first case $Q=M$ and in the second $Q=M^E$. It remains to show
that the remaining cases are impossible. If $Q$ contains $XB$ then $Q=\langle XB,YC,ZA\rangle M^2$ and
$[Q,Q]=\langle X^4,Y^4,Z^4\rangle$. If $Q$ contains $XC$, then $Q=\langle XC,YA,ZB\rangle M^2$ and $[Q,Q]=
\langle X^2Z^2, Y^2X^2, Y^2Z^2, X^4,Y^4,Z^4\rangle$. If $Q$ contains $XAB$ then $Q^E$ is
also in the $\mathrm{Aut}(G)$-orbit of $M$ and contains $XB$ as in the above case. 
If $Q$ contains $XAC$ then $Q^E$ is
also in the $\mathrm{Aut}(G)$-orbit of $M$ and contains $XC$ as in the above above. If $Q$ contains $XABC$ then
it also contains $X,Y,Z,ABC$ and is at least 4-dimensional. If $Q$ contains $XBC$ then $Q^E$ is
also in the $\mathrm{Aut}(G)$-orbit of $M$ and contains $XABC$ as in the above case. All of these cases are impossible.

This proves that the $\mathrm{Aut}(G)$-orbit of $M$ consists of $M$ and $M^E$ only.

\section{Finding $\mathrm{Aut}(G)$ and $\mathrm{Out}(G)$}

We can now prove that $P=\mathrm{Aut}(G)=\mathrm{Inn}(G)\langle E\rangle$ and $\mathrm{Out}(G)\cong C_2$.

Given any $\alpha\in \mathrm{Aut}(G)$,
multiplying it by $E$, we may assume that $\alpha$ preserves $M$. Then, 
the automorphism that $\alpha$ induces on $M$ is inner,
as we have all of $G$ at our disposal. Thus, multiplying $\alpha$ by an inner automorphism, we may assume that $\alpha$
fixes $M$ pointwise. It also preserves $R$, with $R^2=M^2$ also fixed pointwise. Now $A\mapsto A^i B^j C^k$,
so $A^2\mapsto A^2=A^{2i} B^{2j} C^{2k}$, whence $i=1$, $j=0$, $k=0$, and $A$ is fixed. Likewise, $B$ and $C$ are fixed. 

Set $D_0=D^\alpha$. Then $D_0=X^i Y^j Z^k A^a B^b C^d D^\ell$, where $a,b,c\in\Z$, $0\leq i,j,k\leq 1$,
and $\ell\in\{1,2\}$ ($\ell$ cannot be 0 as $U$ is characteristic). 
We must have $A^{D_0}=B$. If $i=1$ then $A^{D_0}=(A^{-1})^{D^\ell}$, which is impossible.
Thus $i=0$ and likewise $j=k=0$. Thus $D_0=A^a B^b C^c D^\ell$. Now $X^{D_0}=Y$, and this forces $a=0$, and
likewise we see that $b=c=0$. But then $\ell=1$.

This proves that $P=\mathrm{Aut}(G)=\mathrm{Inn}(G)\langle E\rangle$. As $E\notin \mathrm{Inn}(G)$ and $E^2$ is conjugation by $ABC$,
it follows that $\mathrm{Out}(G)\cong C_2$.

\section{Completeness of $\mathrm{Aut}(\mathrm{Aut}(M))$}

We identify $G$ with $\mathrm{Inn}(G)$ inside of $\mathrm{Aut}(G)$, so that $\mathrm{Aut}(G)$ is the group generated by
the elements $X,Y,Z,A,B,C,D,E$ subject to the defining relations (\ref{defG1}), (\ref{defG2}), as well as
$$
A^E= A, B^E=B, C^E=C, D^E=D, X^E=XA, Y^E=YB, Z^E=ZC.
$$

We claim that $G$ a characteristic subgroup of $P=\mathrm{Aut}(G)$.  Indeed, we have 
$$
[X,E]=A,
$$
and it follows that
$$
[P,P]=\langle A,B,C, XY, YZ, ZX\rangle,
$$
with $P/[P,P]\cong C_2\times C_2\times C_3$. Thus, there are exactly 3 subgroups of $P$ of index 2, namely
$$G=[P,P]\langle X,D\rangle, G_1=[P,P]\langle E,D\rangle, G_2=[P,P]\langle EX,D\rangle.$$
Now
$$
[G,G]=[G_1,G_1]=[G_2,G_2]=\langle AB,BC,AC,XY, YZ, ZX\rangle,
$$
with factors isomorphic to $C_2\times C_2\times C_3$, so each of $G,G_1,G_3$ has a unique subgroup of index 3, respectively
equal to
$$
U=\langle A,B,C,X,Y,Z\rangle, U_1=\langle A,B,C,XY,YZ,ZX,E\rangle, U_2=\langle A,B,C,XY,YZ,ZX,XE\rangle.
$$
Here
$$
[U,U]=M^2, [U_1,U_1]=[U_2,U_2]=\langle AB,BC,AC\rangle M^2.
$$
But 
$$
[P,[P,P]]=[G,G],
$$
so
$$
[P,[P,P]]^2=[G,G]^2=M^2
$$
is a characteristic subgroup of $P$. If $G$ could be mapped via an automorphism of $P$ to $G_1$ then $U$ would be mapped to $U_1$,
and hence $[U,U]=M^2$ to $[U_1,U_1]=\langle AB,BC,AC\rangle M^2$, which is impossible since $M^2$ is mapped to itself.
An analogous argument shows that $G$ cannot be mapped to $G_2$. Thus $G$ is characteristic subgroup of $P$.

We now appeal to a well-known result of Burnside (see \cite[Chapter II]{Z}) to the effect that if $H$ is a centerless group
such that $\mathrm{Inn}(H)$ is a characteristic subgroup of $\mathrm{Aut}(H)$, then $\mathrm{Aut}(H)$
is a complete group. Applying this to $H=G=\mathrm{Aut}(M)$, we deduce that $\mathrm{Aut}(\mathrm{Aut}(M))$ is a complete group.


\end{document}